\documentclass[a4paper 11pt]{article}
\usepackage{latexsym,amsmath, amsthm, amssymb}
\catcode`@=11 \long\def\@makefntext#1{\noindent #1}
\newskip\tabcentering \tabcentering=1000pt plus 1000pt minus 1000pt
\def\MCH#1#2{\setbox0=\hbox{\raise#1\hbox{#2}}\smash{\box0}}% move char

\def\@evenfoot{}\def\@oddfoot{}

\def\@evenhead{\hbox to\textwidth{\footnotesize\rm\thepage \hfill
{\it }}}

\def\@oddhead{\hbox to \textwidth{\footnotesize{\it
} \hfill\thepage}}

\def\proof{\vspace{2mm}\noindent{\it Proof}\quad}

\newtheorem{thm}{Theorem}[section]
\newtheorem{definition}{Dfinition}[section]

\newtheorem{pro}{Proposition}[section]

\def\bc{\begin{center}}
\def\ec{\end{center}}

\def\hang{\hangindent\parindent}
\def\textindent#1{\indent\llap{\qquad #1\ \ \enspace}\ignorespaces}
\def\ref{\par\hang\textindent}

\begin{document}
 \abovedisplayskip=8pt plus 1pt minus 1pt
\belowdisplayskip=8pt plus 1pt minus 1pt
%-------------------  First Head  -----------------------------------------
\thispagestyle{empty} \vspace*{-3.0truecm} \noindent
\parbox[t]{6truecm}{\footnotesize\baselineskip=11pt\noindent  {} %Acta Mathematica
%%Sinica, English Series\\
%%1999, Jan., Vol.15, No.1, p. 1--11\\
%%Http://www.ActaMath.com\\
%DOI:
 } \hfill
%%\parbox[c]{6truecm}{\vbox{\hsize 3.6576 true cm %
%%  \vskip 3.8 true cm %1.8373
%%  \relax\hbox to0.4\hsize{\hbox to0pt{\special{BMF=actmark.BMF}}\hss}\hss}}
%\hbox to\textwidth{\vbox{\footnotesize\baselineskip=11pt\noindent
%Acta Mathematica Sinica, English Series\hfill LOGO\\
%1999, Jan., Vol.15, No.1, pp. 1--11\hfill \copyright Spring-Verlag 1999}
%\vbox{\hsize3.6576 true cm
%  \vskip1.8373 true cm
%  \relax\hbox to\hsize{\hbox to0pt{\special{BMF=ACTMARK.BMF}}\hss}\hss}}
%===================Text=============================================
\vspace{1 true cm}

\bc{\large\bf  Inf-convolution of $g_\Gamma$-solution and its
applications
  }
%\footnotetext{\footnotesize Received February 24, 1998, Revised
%September 1, 1998, Accepted September 9, 1996}}
\footnotetext{\footnotesize + Corresponding author.\\
} \ec

\vspace*{0.1 true cm}
\bc{\bf Yuanyuan SUI, Helin Wu$^{+}$ \\
{\small\it School of Mathematics, Shandong University, Jinan 250100, China\\
\small\it \quad E-mail: wuhewlin@gmail.com}}\ec \vspace*{3 true mm}

\begin{abstract}
A risk-neutral method is always used to price  and hedge  claims in
complete market, but another wildly used method  in more general
case is based on utility maximization or risk minimization. All
kinds of risk measure have been used in literature. In this paper,
We use a kind of risk measure induced by $g_\Gamma$-solution  or the
minimal solution of a Constrained  Backward Stochastic Differential
Equation (CBSDE) when constraints in investment comes to our
consideration. We adopt the inf-convolution of convex risk measures
to solve some optimization problem.  A dynamic version risk measures
defined through $g_\Gamma$-solution will be get. Just like the case
without constraint, the inf-convolution  of two minimal solutions of
CBSDE with two different coefficients is equivalent to  that of
CBSDE with the inf-convolution of their coefficients. In this case,
it is also possible to characterize the optimal risk transfer.

\end{abstract}

\section {Introduction}
The theory of Backward Stochastic Differential Equation (shortly
BSDE)  and risk measure are two wonderful tools to price and hedge
claims in financial market. Useful reference about these can be
found in Pardox and Peng [9] and Artzner et al. [2] and Delbaen [3];
F\"{o}llmer and Schied [5], [6], [7] and  Frittelli and Rosazza [8],
[9]. Unsurprisingly, one may wonder  if there is some relationship
between them, fortunately, Rossazza [4] has done this work, that is
some kind of useful risk measure can be induced by g-expectation.

In a complete market, a kind  of risk-neutral method is always used
to price and hedge claims via equivalent martingale measure.
However, when the market is incomplete or more generally  when some
constraints were put on wealth and portfolio process, one need to
use super-hedging strategy to get upper price. In this paper, we
define a risk measure via $g_\Gamma$-solution, which is a newly
notation given by the author in Peng and Xu [14].  Interestingly, we
prove such risk measure satisfies the important Fatou-property.

The risk measure induced by $g_\Gamma$-solution is different from
the market modified risk measure used in Pauline Barrieu.,  Nicole
El Karoui [11], [12]. In their paper, a market modified risk measure
was defined as a inf-convolution of some risk measure and the risk
measure generated by some convex set which usually can be viewed as
some constraints in hedging problem. To make the  risk measure
generated by  some set  be well defined, one always ask the set to
satisfy some additional conditions. A convenience  to use the risk
measure induced by $g_\Gamma$-solution is that we need not such
conditions any more.

This paper is organized as follows: In section 2, we state the
 framework in Peng[13] and some propositions about
 $g_\Gamma$-solution.  Under  some mild assumptions, $g_\Gamma$-solution
 is well defined on  $L^\infty(\mathcal{F})$, the space of (P)-essentially bounded variables on some
probability  space $(\Omega,\mathcal{F},P)$. Some results about the
risk measure induced by such solution and some applications of it
are given in section 3.

\section {BSDE  and $g_\Gamma$-solution of CBSDE}

\noindent Given a probability space $(\Omega,\mathcal{F},P)$ and
$R^d$-valued Brownian motion \mbox{$W(t)$}, we consider a sequence
$\{(\mathcal{F}_t);t\in[0,T]\}$ of filtrations generated by Brownian
motion $W(t)$  and augmented by P-null sets. $\mathcal{P}$  is the
$\sigma$-field of predictable sets of  $\Omega\times[0,T]$. We use
$L_T^2(R^d)$ to denote the space of all $F_T$-measurable random
variables $\xi:\Omega\rightarrow R^d$ for which
$$\parallel \xi\parallel^2=E[|\xi|^2]<+\infty.$$
and use $H_T^2(R^d)$ to denote  the space of predictable process
$\varphi:\Omega\times[0,T]\rightarrow R^d$ for which
$$\parallel \varphi\parallel^2=E[\int_0^T|\varphi|^2]<+\infty.$$
The backward stochastic differential equation (shortly BSDE ) driven
by $g(t,y,z)$ is given by
$$-dy_t=g(t,y_t,z_t)dt-z^*_tdW(t) \eqno(2.1)$$
 where $y_t\in R$ and $W(t)\in R^d$.
Suppose that $\xi\in L_T^2(R)$ and $g$ satisfies
$$|g(\omega,t,y_1,z_1)-g(\omega,t,y_2,z_2)|\leq M(|y_1-y_2|+|z_1-z_2|), \ \ \forall
(y_1,z_1),(y_2,z_2)\eqno(A1)$$ for some $M>0$  and
$$g(\cdot,0,0)\in H_T^2(R) \eqno(A2)$$

Pardoux and Peng [10] proved the existence  of adapted solution
$(y(t),z(t))$ of such BSDE. We call $(g,\xi)$ standard parameters
for the BSDE.

The following definitions is necessary to help us go on with our
study.

\begin{definition}
(super-solution) A super-solution of a BSDE associated with the
standard parameters $(g,\xi)$ is a vector process $(y_t,z_t,C_t)$
satisfying
$$-dy_t=g(t,y_t,z_t)dt+dC_t-z^*_tdW(t),\quad y_T=\xi,\eqno(2.2)$$
or being equivalent to
$$y_t=\xi+\int_t^Tg(s,y_s,z_s)ds-\int_t^Tz^*_sdW_s+\int_t^TdC_s, \eqno(2.2')$$
where $(C_t,t\in[0,T])$ is an increasing, adapted, right-continuous
process with $C_0=0$ and $z_t^*$ is the transpose  of $z_t$. When
$C_t\equiv 0$, we call $(y_t, z_t)$ a g-solution.
\end{definition}
\noindent Constraints like
$$(y(t),z(t))\in \Gamma\eqno(C)$$
where $\Gamma=\{(y,z)|\phi(y,z)=0\}\subset R\times R^d$ and
$\phi(y,z): R\times R^d\rightarrow R^+$ is always  considered  in
this paper. In such case, we  give the following definition,
\begin{definition}( $g_\Gamma$-solution or the minimal solution ) A g-supersolution $(y_t, z_t, C_t)$ is said to be the
 the minimal solution, given $y_T=\xi$,
subjected to the constraint $(C)$ if for any other g-supersolution
$(y'_t, z'_t, C'_t)$ satisfying $(C)$  with $y'_T=\xi$, we have
$y_t\leq y'_t $ a.e., a.s.. The  minimal solution is denoted by
$\mathcal{E}_t^{g,\phi}(\xi)$ and for convenience called as
$g_\Gamma$-solution.
\end{definition}
For any $\xi\in L^2_T(R)$, we denote  $\mathcal{H}^{\phi}(\xi)$ as
the set of g-supersolutions $(y_t,z_t,C_t)$ subjecting to $(C)$ with
$y_T=\xi$.  When $\mathcal{H}^{\phi}(\xi)$ is not empty, Peng [13]
proved that $g_\Gamma$-solution exists.

The convexity of $\mathcal{E}_t^{g,\phi}(\xi)$ can be easily deduced
from the same proposition of solution of BSDE with convex generator
function.

\begin{pro}
Let  $\phi(t,y,z)$ be a function: $[0,T]\times R\times
R^d\rightarrow R^+$ and $g(t,y,z)$ be a  function: $[0,T]\times
R\times R^d\rightarrow R$. Suppose $\phi(t,y,z)$ and $g(t,y,z)$ are
both convex in $(y,z)$ and satisfy (A1) and (A2), then
$$\mathcal{E}_t^{g,\phi}(a\xi+(1-a)\eta)\leq a\mathcal{E}_t^{g,\phi}(\xi)+(1-a)\mathcal{E}_t^{g,\phi}(\eta)\quad \forall t\in[0,T]$$
holds for any $\xi,\eta\in L_T^2(R)$ and $a\in[0,1]$.
\end{pro}

\proof According to Peng [13], the solutions  $y_t^m(\xi)$ of
$$y_t^m(\xi)=\xi + \int_t^Tg(y_s^m(\xi),z_s^m,s)ds+A_T^m-A_t^m-\int_t^Tz_s^mdW_s.$$
is an increasing sequence and converges  to
$\mathcal{E}_t^{g,\phi}(\xi)$, where $$A_t^m: =
m\int_0^t\phi(y_s^m,z_s^m,s)ds.$$ For any fixed $m$, by the
convexity of $g$ and $\phi$, $y_t^m(\xi)$ is a convex in $\xi$, that
is
$$y_t^m(a\xi+(1-a)\eta)\leq ay_t^m(\xi)+(1-a)y_t^m(\eta),$$
taking limit as $m\rightarrow\infty$, we get the required result.
\hspace*{\fill}$\Box$

By the same method of penalization, we can get the comparison
theorem of $\mathcal{E}_t^{g,\phi}(\xi)$ .
\begin{pro}
Under the same assumptions as above proposition, we have
$$\mathcal{E}_t^{g_1,\phi}(\eta )\geq
\mathcal{E}_t^{g_2,\phi}(\xi), \quad\forall t\in[0,T] \quad P-a.s.$$
for any $\xi, \eta\in L_T^2(R)$ when $P(\eta\geq\xi)=1$ and $g_1\geq
g_2$.

\end{pro}

\section {Risk measure via $g_\Gamma$-solution  and its applications}
In this section, we study convex risk measure induced by
$g_\Gamma$-solution. First we give the concept of convex risk
measure  which can be got from many papers such as F\"{o}llmer and
Schied [5].
\begin{definition}
Let $L^\infty(P)$ be  the space of (P)-essentially bounded functions
on some probability  space $(\Omega,\mathcal{F},P)$. A functional
$\rho: L^\infty(P)\longrightarrow R$ is a (moneytary) convex risk
measure if, for any $\xi$ and $\eta$ in  $L^\infty(P)$, it satisfies
the following properties:

a)Convexity: $\forall \lambda\in[0,1] \qquad
\rho(\lambda\xi+(1-\lambda)\eta)\leq\lambda\rho(\xi)+(1-\lambda)\rho(\eta);$

b) Monotonicity: $\xi\leq\eta \quad a.s(P)\Rightarrow
\rho(\xi)\geq\rho(\eta);$

c)Translation invariance: $\forall m\in
R\quad\rho(\xi+m)=\rho(\xi)-m.$

A convex risk measure $\rho$ is coherent if it satisfies also:

d)Homogeneity:$\quad\forall \lambda\in R^+ \qquad
\rho(\lambda\xi)=\lambda\rho(\xi).$

\end{definition}

 In order to generate a convex risk measure by
$g\Gamma$-solution, we need some additional assumptions such as
$$g \,\,\text{is independent of}\,\, y \,\, \text{and }\,\, g(\cdot,0)=0 \eqno(A3)$$

When $g$ satisfying conditions $A(i), i=1,2,3$, just as Rosazza [4]
noted, some useful risk measure can be generated by $g$-expectation.

First, we prove a result that $g_\Gamma$-solution can be well
defined on the space  $L^\infty(\mathcal{F}_T)$ of (P)-essentially
bounded functions on some probability  space
$(\Omega,\mathcal{F_T},P)$.

\begin{pro}
Suppose that $g$ and $\phi$ satisfy  assumptions $A(i), i=1,2,3$,
then $\mathcal{E}_t^{g,\phi}(\cdot)$ is  well defined on
$L^\infty(\mathcal{F}_T)$

\end{pro}
\proof Since $g$ is independent of $y$ and $g(t,0)=0$,
$\phi(t,0)=0$, then for any  fixed $C_0>0, \mu>0$, we have
$$g(t,y,0)\leq C_0+\mu|y|, \quad (y,0)\in \Gamma_t,\quad \forall
y\geq C_0.$$

By  Peng and Xu [14], the $g_\Gamma$-solution with terminal
condition $y_T=\xi$ exists for any $\xi\in
L^2_{+,\infty}(\mathcal{F}_T)$,  where

$$ L^2_{+,\infty}(\mathcal{F}_T):=\{\xi\in L^2(\mathcal{F}_T),
\xi^+\in L^\infty(\mathcal{F}_T)\}.$$

It is obvious $L^\infty(\mathcal{F}_T)\subset
L^2_{+,\infty}(\mathcal{F}_T),$  thus $\mathcal{E}_t^{g,\phi}(\xi)$
exists for any $L^\infty(\mathcal{F}_T)$. \hspace*{\fill}$\Box$

We first consider the case $t=0$, then
$\rho(\xi)=\mathcal{E}_0^{g,\phi}(-\xi )$ generated a static convex
risk measure when both $g$ and $\phi$ are convex functions
satisfying assumptions $A(i),i=1,2,3$. Furthermore, we can prove
$\rho$ satisfies the important Fatou property.

\begin{thm}
When both $g$ and $\phi$   satisfy assumptions $A(1)$ and $A(2)$,
then $\mathcal{E}_0^{g,\phi}(\xi )$ is continuous from below, etc,
when $\{\xi_n\in L^\infty(\mathcal{F}_T), \ n=1,2,\cdots\}$ is an
increasing sequence  comes from $L^\infty(\mathcal{F}_T)$ and
converges almost surely to $\xi\in L^\infty(\mathcal{F}_T)$, then
$$\lim_{n\rightarrow \infty}\mathcal{E}_0^{g,\phi}(\xi_n)= \mathcal{E}_0^{g,\phi}(\xi).$$
\end{thm} \proof Taking  $y_t^m(\xi)$  as  in proposition 2.1.
  By proposition 2.2,
$\{\mathcal{E}_t^{g,\phi}(\xi_n),n=1,2,\cdots\}$ is an increasing
sequence. We denote its limit at $t=0$ as  $a$,  then $a\leq
\mathcal{E}_0^{g,\phi}(\xi)$. Since $\xi_n$ converges almost surely
increasingly to $\xi\in L^\infty(\mathcal{F}_T)$, by dominated
convergence theorem, it also converges strongly  in $L_T^2(P)$, then
by the continuous dependence property of g-supersolution, the limit
of $\{y_0^m(\xi_n)\}_{n=1}^{\infty}$ is $y_0^m(\xi)$ for any fixed
$m$.

We want to show  that $a=\mathcal{E}_0^{g,\phi}(\xi)$. If on the
contrary on has  $a<\mathcal{E}_0^{g,\phi}(\xi)$, then there is some
$\delta>0$ such that
$\mathcal{E}_0^{g,\phi}(\xi)-\mathcal{E}_0^{g,\phi}(\xi_n)>\delta$
for any $n$. On the other hand, for any $\epsilon>0$, $0\leq
\mathcal{E}_0^{g,\phi}(\xi)-y_0^m(\xi)\leq\epsilon$ holds for some
larger $m_0$. Fixing  $m_0$, $\epsilon$, there is some $n_0$ which
depends on $m_0$ and $\epsilon$ such that $0\leq y_0^{m_0}(\xi)-
y_0^{m_0}(\xi_{n_0})\leq\epsilon$, so
$\mathcal{E}_0^{g,\phi}(\xi)-y_0^{m_0}(\xi_{n_0})\leq 2\epsilon$,
but we have $\mathcal{E}_0^{g,\phi}(\xi)-y_0^{m_0}(\xi_{n_0})\geq
\mathcal{E}_0^{g,\phi}(\xi)-\mathcal{E}_0^{g,\phi}(\xi_{n_0})>\delta$,
this is impossible for $\epsilon<\frac{\delta}{2}$.
\hspace*{\fill}$\Box$

Thanks to this property and the work done by F\"{o}llmer, H., Schied
[6], [7], the convex risk measure can be represented by a family of
probabilities which are absolutely continuous with $P$.

We then go to some applications of $g_\Gamma$-solution. Here we use
some notations in Pauline Barrieu.,  Nicole El Karoui [11]. Let
$\eta\in L_T^\infty(P)$, $\rho(\xi)=\mathcal{E}_0^{g,\phi}(-\xi )$
be a convex risk measure when both $g$ and $\phi$ are convex, our
first problem is  a minimizing  problem by inf-convolution. More
explicitly, suppose two agents who have convex risk measure
generated by $\rho_i(\xi)=\mathcal{E}_0^{g_i,\phi_i}(-\xi ),i=1,2$
respectively, we want to find an optimal value in $L_T^\infty(P)$ to
attain
$$\inf_{\xi\in L^\infty(\mathcal{F}_T)}\{\rho_1(\eta-\xi)+\rho_2(\xi)\}.\eqno(3.1)$$

We first consider two simple cases.
\begin{thm}
If both $g$ and $\phi$  satisfy assumptions $A(i),i=1,2,3$ and
$$h(z_1+z_2)\leq h(z_1)+h(z_2), \forall z_1,z_2 $$
holds for $h=g, \phi$, then $\xi=0$ is a optimal value for problem
(3.1) when $g_i=g,\phi_i=\phi,i=1,2$.
\end{thm}
\proof Suppose that $(y(t)=\mathcal{E}_t^{g,\phi}(\xi
-\eta),z(t),C(t))$ and  $(\tilde{y}(t)=\mathcal{E}_t^{g,\phi}(-\xi
),\tilde{z}(t),\tilde{C}(t))$ are $g_\Gamma$-solutions with terminal
value $\xi-\eta$ and $-\xi$ respectively, that is
$$y_t=\xi-\eta+\int_t^Tg(s,z_s)ds-\int_t^Tz^*_sdW_s+\int_t^TdC_s, \eqno(3.2)$$
$$\tilde{y}_t=-\xi+\int_t^Tg(s,\tilde{z}_s)ds-\int_t^T\tilde{z}^*_sdW_s+\int_t^Td\tilde{C}_s. \eqno(3.3)$$
Add (3.2) and (3.3) together, we have
$$\mathcal{E}_t^{g,\phi}(\xi-\eta)+\mathcal{E}_t^{g,\phi}(-\xi)\geq
-\eta+\int_t^T(g(s,z_s)+g(s,\tilde{z}_s))ds-\int_0^T(z^*_s+\tilde{z}^*_s)dW_s+\int_t^Td(C_s+\tilde{C}_s).\eqno(3.4)$$

By the assumption, we  have furthermore that
$$y(t)+\tilde{y}(t)\geq
\bar{y}(t):=-\eta+\int_t^Tg(s,z_s+\tilde{z}_s)ds-\int_0^T(z^*_s+\tilde{z}^*_s)dW_s+\int_t^Td(C_s+\tilde{C}_s).\eqno(3.5)$$
and $0\leq\phi(z_s+\tilde{z}_s)\leq \phi(z_s)+\phi(\tilde{z}_s)=0.$

This means that $(\bar{y}(t),z(t)+\tilde{z}(t),C(t)+\tilde{C}(t))$
is a super-solution  with terminal value $-\eta$. By (3.5) and the
definition of $g_\Gamma$-solution, we have
$$\mathcal{E}_t^{g,\phi}(\xi-\eta)+\mathcal{E}_t^{g,\phi}(-\xi)\geq
\mathcal{E}_t^{g,\phi}(-\eta).$$

Take $t=0$, we have
$$\rho(\eta-\xi)+\rho(\xi)\geq\rho(\eta),\quad \forall \xi\in L_T^\infty(P).$$
This means $\xi=0$ is a optimal value for problem
(3.1).\hspace*{\fill}$\Box$

The result of above tells us that if two agents having risk measure
induced by same coefficients, then one rational  way of them to
transfer risk is doing nothing.

We then go to consider another interesting case concerning  a useful
operator of risk measure. For any $\lambda>0$, which always be
considered as the risk tolerance coefficient, we can define the
dilatation of convex risk measure $\rho(\xi)$  as
$\rho_\lambda=\lambda\rho(\xi/\lambda)$. Our first result is that
under some mild assumptions, the dilatation risk measure of
$g_\Gamma$-solution coincides with the minimal solution of the
dilatation of coefficients.
\begin{thm}
Suppose $g$ and $\phi$ satisfy the assumptions $A(i),i=1,2,3$,
 $\phi(\lambda z)=\lambda\phi(z)$ holds for any $0<\lambda$. Let $\rho(\xi)=
\mathcal{E}_0^{g,\phi}(-\xi)$, $g_\lambda(z)=\lambda g(z/\lambda)$,
then we have
$$\lambda\rho(\xi/\lambda)=
\mathcal{E}_0^{g_\lambda,\phi}(-\xi)$$
\end{thm}
 \proof Suppose that $(y(t),z(t),C(t))$ is the $g_\Gamma$-solution
 with terminal value $\xi/\lambda$,
 $$\mathcal{E}_t^{g,\phi}(\xi/\lambda)=y_t=\xi/\lambda+\int_t^Tg(s,z_s)ds-\int_t^Tz^*_sdW_s+\int_t^TdC_s.\eqno(3.6)$$
  then $$\lambda\mathcal{E}_t^{g,\phi}(\xi/\lambda)=\lambda y_t=\xi+\int_t^T\lambda g(s,z_s)ds-\int_t^T\lambda z^*_sdW_s+\int_t^Td\lambda C_s\eqno(3.7)$$
At the same time we suppose that
$(\tilde{y}(t),\tilde{z}(t),\tilde{C}(t))$ is the minimal solution
with coefficient $g_\lambda=\lambda g(z/\lambda)$ and terminal value
$\xi$ satisfying constraint (C),
$$\mathcal{E}_t^{g_\lambda,\phi}(\xi)=\tilde{y}_t=\xi+\int_t^Tg_\lambda(s,\tilde{z}_s)ds-\int_t^T\tilde{z}^*_sdW_s+\int_t^Td\tilde{C}_s.\eqno(3.8)$$
By (3.7), we can see that $(\lambda y_t, \lambda z_t, \lambda C_t)$
is a $g_\lambda$-supersolution with terminal value $\xi$ satisfying
constraint $(C)$, thus we have
$$\lambda \mathcal{E}_t^{g,\phi}(\xi)\geq
\mathcal{E}_t^{g_\lambda,\phi}(\xi) \quad  a.e\quad a.s.\eqno(3.9)$$
Similarly, by (3.8),
$(\tilde{y}(t)/\lambda,\tilde{z}(t)/\lambda,\tilde{C}(t)/\lambda)$
is a $g$-supersolution with terminal  value $\xi$ satisfying
constraint (C), thus we have
$$\mathcal{E}_t^{g,\phi}(\xi)\leq
\mathcal{E}_t^{g_\lambda,\phi}(\xi)/\lambda \quad  a.e\quad
a.s.\eqno(3.10)$$

Put (3.9) and (3.10) together, we get

$$\lambda \mathcal{E}_t^{g,\phi}(\xi)=
\mathcal{E}_t^{g_\lambda,\phi}(\xi) \quad  a.e\quad a.s.$$

Specially $$\lambda\rho(\xi/\lambda)=
\mathcal{E}_0^{g_\lambda,\phi}(-\xi)$$ holds.\hspace*{\fill}$\Box$

Thanks to this result and the wonderful result in Pauline Barrieu.,
Nicole El Karoui [11], we have the following result.
\begin{thm}
Suppose $g$ and $\phi$ satisfy the assumptions $A(i),i=1,2,3$, two
agents  have risk measure with different risk tolerance coefficient
$g_\lambda$ and $g_\gamma$ respectively, then one optimal value of
problem (3.1) is
$$\xi=\frac{\gamma}{\gamma+\lambda}\eta.$$
\end{thm}

When one consider the optimal problem (3.1) with general
coefficients $g_i,i=1,2$, we need more concepts.
\begin{definition}
 Let $X$ be a Banach space, $X^*$ is its dual space and  $\varphi :X\rightarrow R$ is a  convex functional. For any $\xi \in X$,
 define
 $$\partial\varphi(\xi)\triangleq \{f\in X^*,f(\eta)\leq\varphi(\xi+\eta)-\varphi(\xi),\forall \eta\in X\}$$
 as the subdifferential  of $\varphi$ at $\xi$,
 every member  of $\partial\varphi(\xi)$ is called a subgradient of subdifferential  of $\varphi$ at $\xi$.
\end{definition}
The following result is  basic  in convex analysis.

\begin{pro}
Suppose  $\varphi$ is a continuous   convex functional on $X$, then
for any $\xi\in X$,  $\partial\varphi(\xi)$ is  not empty.
\end{pro}

 \proof In the product space $X\times R$, let
 $D\triangleq \{(\xi,t)|\varphi(\xi)\leq t\}$ be the
 upper semi-graph of $\varphi$. For any fixed point $\xi_0\in X$,
 since $\varphi(\cdot)$ is continuous at $\xi_0$, $(\xi_0,
 \varphi(\xi_0)+1)$ is a interior point of $D$.
 Note that
 $$\{(\xi_0,\varphi(\xi_0))\}\bigcap \mathring{D}={\O},$$
 then by separating  theorem of convex sets in Banach space,
 there is some  no zero point $(g,a)\in X^*\times R$ such that
 $$g(\xi_0)+ a\varphi(\xi_0)\leqslant g(x)+ at \quad \forall (x,t)\in
 D.$$
It is not hard to check  that $a>0$, then if we take $f=-g/a$, then
$f\in\partial\varphi(\xi)$.\hspace*{\fill}$\Box$

The next result gives us a sufficient condition for a convex
functional to be continuous, for its proof, we refer to Aubin[1].
\begin{pro}
 Let $X$ be a Banach space, $\varphi :X\rightarrow R$ be a  convex
 functional. If $\varphi$ is lower semi-continuous on $X$, then it
 is continuous  on $X$.
\end{pro}
 A useful result has been obtained in our previous paper.
 \begin{thm}
Suppose $g$ and $\phi$ satisfy the assumptions $A(i),i=1,2,3$, then
$\mathcal{E}_0^{g,\phi}(\xi)$ is lower semi-continuous on
$L^\infty(\mathcal{F}_T)$.
 \end{thm}
 \proof See Wu [15] for reference.\hspace*{\fill}$\Box$

 We then have a general result when two agents have risk measure
 generated  by general coefficients $g_i,\phi_i,i=1,2.$

 \begin{thm}
Suppose $g_i,\phi_i,i=1,2$ are convex functions satisfying the
assumptions $A(i),i=1,2,3$ and there is some $a,b\in R$ such that
$g_i(t,z)\geq az+b,i=1,2$. If there is some $\xi\in
L^\infty(\mathcal{F}_T)$ and some finite additive measure
$Q\in\partial\hat{\rho}(\eta)\bigcap\partial\rho_1(\eta-\xi)\bigcap\partial\rho_2(\xi)$,
then $\xi$ is optimal for problem (3.1), where
$$\rho_i(\cdot)=\mathcal{E}_0^{g_i,\phi_i}(-\cdot),i=1,2;\quad \hat{\rho}(\cdot)=\inf_{\xi\in
L^\infty(\mathcal{F}_T)}\{\rho_1(\cdot-\xi)+\rho_2(\xi)\}.$$
 \end{thm}
\proof By the assumption that $g_i(z)\geq az+b,i=1,2$, we have that
the inf-convolution $\hat{\rho}$ is well defined on
$L^\infty(\mathcal{F}_T)$. The rest of the proof can be found in
Pauline Barrieu.,  Nicole El Karoui [12].

At last, we state a dynamic version of inf-convolution of
$g_\Gamma$-solution.
\begin{thm}
Suppose $g_i,i=1,2, \phi$ are convex functions satisfying the
assumptions $A(i),i=1,2,3$,
$\phi(t,z_1+z_2)\leq\phi(t,z_1)+\phi(t,z_1),\,\forall z_1,z_2$ and
there is some $a,b\in R$ such that $g_i(t,z)\geq az+b,i=1,2$. The
inf-convolution of $g_1$ and $g_2$ is given by
$$g_3(t,z)=g_1\square g_2(t,z)=\inf_{y}\{g_1(t,z-y)+g_2(t,y)\}.$$
Let  $(\mathcal{E}_t^{g_3,\phi}(\eta), \hat{z}_3(t), \hat{C}_3(t))$
be the $g_\Gamma$-solution with terminal value $\xi\in
L^\infty(\mathcal{F}_T)$ satisfying constraint (C) and $\hat{z}$ be
a measurable process such that
$\hat{z}=\arg\min_y\{g_1(t,\hat{z}_3(t)-y)+g_2(t,y)\} \quad dt\times
dP-a.s.$, then the following results hold:

(1) For any $t\in[0,T]$ and any $\xi\in L^\infty(\mathcal{F}_T)$,
$$\mathcal{E}_t^{g_3,\phi}(\eta)\leq
\mathcal{E}_t^{g_1,\phi}(\eta-\xi)+\mathcal{E}_t^{g_2,\phi}(\xi).$$

(2) If $\phi(t,\hat{z}(t))=0, \,\phi(t,\hat{z}_3(t)-\hat{z}(t))=0$
and

$$\xi^*:=\int_0^Tg_2(s,\hat{z}_s)ds-\int_0^T\hat{z}^*_sdW_s\in
L^\infty(\mathcal{F}_T),$$

then $\xi^*$ is an optimal value for problem (3.1), furthermore, we
have
$$\mathcal{E}_t^{g_3,\phi}(\eta)=
\mathcal{E}_t^{g_1,\phi}\square\mathcal{E}_t^{g_2,\phi}(\eta), \quad
\forall t\in[0,T].$$

\end{thm}
\proof (1) By the same argument of proposition 3.1,
$\mathcal{E}_t^{g_3,\phi}(\eta)$ exists for any $\eta\in
L^\infty(\mathcal{F}_T)$.

Suppose that $(y_i(t),z_i(t),C_i(t)),i=1,2$ is the minimal solution
 with terminal value $\eta-\xi$ and $\xi$ for CBSDE with
 coefficients $g_i$ satisfying constraint (C), that is
 $$\mathcal{E}_t^{g_1,\phi}(\eta-\xi)=y_1(t)=\eta-\xi+\int_t^Tg(s,z_1(s))ds-\int_t^Tz^*_1(s)dW_s+\int_t^TdC_1(s).\eqno(3.11)$$
$$\mathcal{E}_t^{g_2,\phi}(\xi)=y_2(t)=\xi+\int_t^Tg(s,z_2(s))ds-\int_t^Tz^*_2(s)dW_s+\int_t^TdC_2(s).\eqno(3.12)$$

Put (3.11) and (3.12) together, by the comparison property of
proposition 2.2, we have
$$\mathcal{E}_t^{g_1,\phi}(\eta-\xi)+\mathcal{E}_t^{g_2,\phi}(\xi)\geq
y_3(t)=\eta+\int_t^Tg_3(s,z_3(s))ds-\int_t^Tz^*_3(s)dW_s+\int_t^TdC_3(s).$$
where $z_3(t)=z_1(t)+z_2(t), C_3(t)=C_1(t)+C_2(t).$

But $(y_3(t),z_3(t),C_3(t))$ is a $g_3$-supersolution satisfying
constraint (C), we have
$$\mathcal{E}_t^{g_3,\phi}(\eta)\leq
\mathcal{E}_t^{g_1,\phi}(\eta-\xi)+\mathcal{E}_t^{g_2,\phi}(\xi).\eqno(3.13)$$

 (2)

Since
$$\mathcal{E}_t^{g_3,\phi}(\eta)=\eta+\int_t^Tg_3(s,\hat{z}_3(s))ds-\int_t^T\hat{z}^*_3(s)dW_s+\int_t^Td\hat{C}_3(s).\eqno(3.13)$$

But
$g_3(t,\hat{z}_3(t))=g_1(t,\hat{z}_3(t)-\hat{z}(t))+g_2(\hat{z}(t)).$
Let
$$\hat{y}(t)=-\int_0^tg_2(s,\hat{z}(s))ds+\int_0^t\hat{z}^*(s)dW_s,$$
that is
$$\hat{y}(t)=\xi^*+\int_t^Tg_2(s,\hat{z}(s))ds-\int_t^T\hat{z}^*(s)dW_s$$
and it is obvious  that
$\mathcal{E}_t^{g_2,\phi}(\xi^*)=\hat{y}(t).$ By (3.14),
$(\mathcal{E}_t^{g_3,\phi}(\eta)-\mathcal{E}_t^{g_2,\phi}(\xi^*),\hat{z}_3(t)-\hat{z}(t),\hat{C}_3(t))$
is a $g_1$-supersolution with terminal value $\eta-\xi^*$ satisfying
constraint (C), so
$$\mathcal{E}_t^{g_3,\phi}(\eta)-\mathcal{E}_t^{g_2,\phi}(\xi^*)\geq\mathcal{E}_t^{g_1,\phi}(\eta-\xi^*).\eqno(3.15)$$
By (3.13) and (3.15), we get

$$\mathcal{E}_t^{g_3,\phi}(\eta)=\mathcal{E}_t^{g_1,\phi}(\eta-\xi^*)+\mathcal{E}_t^{g_2,\phi}(\xi^*)=\mathcal{E}_t^{g_1,\phi}\square\mathcal{E}_t^{g_2,\phi}(\eta).$$


\begin{thebibliography}{}
\bibitem{1}J.P.Aubin: Optima and Equilibria. An Introduction to
Nonlinear Analysisi. Springer, Berlin, (1998)
\bibitem{2} P.Artzner,  F.Delbaen, J.M.Eber and   D.Heath: Coherent
measures of risk. Mathematical Finance 4, (1999) 203-228
\bibitem{3}F.Delbaen: Coherent risk measures on general probability space. Advance in Finance and Stochastics, springer-verlag, (2002) 1-37

\bibitem{4}Emanuela Rosazza Gianin.: Risk measures via
g-expectations. Insurance Mathematics and Economics 39 (2006) 19¨C34
\bibitem{5}H.F\"{o}llmer and A.Schied: Convex measure of risk and trading constraints.Finance and stochastics 6(4), (2002) 429-447
\bibitem{6}H.F\"{o}llmer and A.Schied: Robust preferences and convex measures of risk. In: Sandmann, K., Sch\"{o}nbucher, P.J. (Eds.), Advances in Finance
and Stochastics. Springer-Verlag, (2002) 39-56
\bibitem{7}H.F\"{o}llmer and A.Schied: Stochastic Finance. An Introduction in Discrete Time. De Gruyter, Berlin, New York.
\bibitem{8} M.Frittelli and E.Rosazza Gianin : Putting order in risk measures. Journal of Banking and Finance 26 (7), (2002) 1473-1486
\bibitem{9} M.Frittelli and E.Rosazza Gianin: Dynamic convex risk measures. In: Szeg\"{o},G.(Ed.), Risk Measures for the 21st Century. J. Wiley,
\bibitem{10}E.Pardoux and S.G.Peng: Adapted solution of a backward stochastic differential equation. Systems and control letters 14 (1990)  55-62.
\bibitem{11}Pauline Barrieu.,  Nicole El Karoui. Inf-convolution of risk measures and optimal risk
transfer. Finance and Stochastics, 9(2) (2005) 269-298,
\bibitem{12}Pauline Barrieu.,  Nicole El Karoui. Pricing, Hedging and Optimally Designing Derivatives Via Minimization of Risk
Measures. Arxiv preprint arXiv:0708.0948, 2007 - arxiv.org
\bibitem{13}Peng, S. Monotonic limit theorem of BSDE and nonlinear decomposition theorem of Doob-Meyer's type. Probab.Theory Relat.Fields. 113 (1999) 473-499.
\bibitem{14}Peng, S., Xu, M.Y. The smallest g-supermartingale and reflected BSDE with single and double $L^2$ obstacles. Ann.I.H.Poincare-PR 41 (2005) 605-630.
\end{thebibliography}
\end{document}